\documentclass[12pt,a4paper]{article}
\usepackage[T1]{fontenc}
\usepackage{amsmath}
\usepackage{amssymb}
\usepackage{textcomp}
\usepackage[english]{babel}
\usepackage[all,dvips]{xy}
\usepackage[applemac]{inputenc}

\newtheorem{ass}{Assumption}
\newtheorem{thm}{Theorem}[section]
\newtheorem{lem}[thm]{Lemma}
\newtheorem{prop}[thm]{Proposition}
\newtheorem{cor}[thm]{Corollary}

\newcommand\dem{\textbf{Proof: }}
\newcommand\Def{\textbf{Definition: }}

\newcommand\rem{\textbf{Remark: }}

\title{Reducibility of quasiperiodic cocycles under a Brjuno-Rssmann arithmetical condition}
\date{}
\author{C.Chavaudret\\
\small{Universit de Nice Sophia-Antipolis, France}\\
and\\
S.Marmi\\
\small{Scuola Normale Superiore, Pisa, Italy}}

\begin{document}
\maketitle

\textbf{Abstract:} The arithmetics of the frequency and of the rotation number play a fundamental role in the 
study of reducibility of analytic quasi-periodic cocycles which are sufficiently close to a constant.
In this paper we show how to generalize previous works by L.H.Eliasson which deal with the diophantine
case so as to implement a  Brjuno-Rssmann arithmetical condition both on the frequency
and on the rotation number. Our approach adapts the Pschel-Rssmann KAM method, which was previously
used in the problem of linearization of vector fields, to the problem of reducing cocycles. 

\tableofcontents

\section{Introduction}

Quasiperiodic cocycles are the fundamental solutions of quasi-periodic linear systems 

\begin{equation}\label{base}\forall (\theta,t)\in \mathbb{T}^d\times \mathbb{R},\ X'(\theta,t)=A(\theta+t\omega)X(\theta,t)
\end{equation}

\noindent 
where $A$ is a continuous matrix-valued function on a torus $\mathbb{T}^d$ and $\omega$ is a rationally independent vector of some space $\mathbb{R}^d$ (the space of frequencies).
Although the dynamics of such a system can be quite complicated, they are easily studied in case the system is  reducible, i.e when there is a map $Z$, continuous on the double torus $2\mathbb{T}^d=\mathbb{R}^d/2\mathbb{Z}^d$, taking its values in the group of invertible matrices and such that

$$\forall \theta \in 2\mathbb{T}^d, \frac{d}{dt}Z(\theta+t\omega)_{\mid t=0}
=A(\theta)Z(\theta)-Z(\theta)B$$

\noindent 
for some matrix $B$ not depending on $\theta$. Since smoothness is an issue, given a class of functions $\mathcal{C}$, we will say that the cocycle is reducible in $\mathcal{C}$ if $Z$ is in 
$\mathcal{C}$.
Here we will focus on the case in which $A$ takes its values in $sl(2,\mathbb{R})$, which is sufficient, for instance, for the study of the one-dimensional quasi-periodic Schrdinger equation.
Moreover, we will consider solutions of \eqref{base} with $A\in C^\omega_r$, the space of functions on $\mathbb{T}^d$ having a holomorphic extension on $\{(z_1,\dots, z_d)\in \mathbb{C}^d,\ \forall j\  \mid \operatorname{Im} z_j\mid <r \}$, whose "weighted norm" $\mid \cdot \mid_r$ converges (see section \ref{defnot}).

%
%

\bigskip
\noindent 
The arithmetics of $\omega$ seem fundamental in the study of reducibility, as well as the arithmetics of the system's rotation number $\rho$ (as it was defined in \cite{JM}; we recall the definition in Section \ref{defnot}). At least in the perturbative case, arithmetical conditions of diophantine type have long been used to obtain reducibility, which can be seen as the convergence of a certain sequence of analytic functions: a diophantine condition can be used to control small divisors and make sure that the sequence converges. This was achieved, in particular, by Eliasson in \cite{El92}:

\begin{thm}[Eliasson] Let $r>0, V\in C^\omega_r(\mathbb{T}^d,\mathbb{R})$. Suppose $\omega$ is a diophantine vector. There exists $\epsilon_0$ depending only on $r,\omega$ such that if 
$\sup_{\mid \operatorname{Im}\theta\mid<r} \mid V(\theta)-\hat{V}(0)\mid\leq \epsilon_0$, then the cocycle which is solution of 

\begin{equation}\label{schr}\frac{d}{dt}X(t,\theta)=\left(\begin{array}{cc}
0 & V(\theta+t\omega) - E\\
1 & 0\\
\end{array}\right) X(t,\theta)
\end{equation}

\noindent is reducible for all $E$ for which the rotation number is rational or diophantine with respect to $\omega$.

\end{thm}


\noindent 
In this article, we will give a reducibility result for analytic cocycles under a weaker arithmetical condition than the diophantine one. In order to obtain an analytic reducibility result, we will have to pick a frequency and a rotation number with good approximation properties, in the sense of Rssmann (\cite{R}): $\omega$ will have to satisfy a strong irrationality condition controlled by an approximation function G, namely

$$\forall m\in\mathbb{Z}^d\setminus\{0\},
\mid \langle m,\omega\rangle\mid\geq \frac{\kappa}{G(m)}$$

\noindent 
for some positive $\kappa$ (section \ref{defnot}), and $\rho$ will have to satisfy a further arithmetical condition: its approximations by means of linear combinations of the frequencies are controlled by an approximation function $g$, i.e

$$\forall m\in\mathbb{Z}^d\setminus\{0\},
\mid\rho - \langle m,\omega\rangle\mid\geq \frac{\kappa'}{g(m)}$$

\noindent (we will say for short that $\rho$ has $g$ as approximation function with respect to $\omega$ with constant $\kappa'$) 
with $g,G$ satisfying some extra assumptions. 

\bigskip
\noindent We will be particularly interested in the case of Brjuno frequency, i.e when 

\begin{equation}\label{brjuno-russmann}\int_1^\infty \frac{\log G(t)}{t^2}dt<\infty
\end{equation}

\noindent and of $\frac{1}{2}$-Brjuno rotation number (with respect to $\omega$), i.e when 

\begin{equation}\label{brjunowrtomega}\int_1^\infty \frac{\log g(t)}{t^{\frac{3}{2}}}dt<\infty
\end{equation}

\noindent 
In dimension $d=2$, Condition \eqref{brjuno-russmann} coincides with the well-known Brjuno condition defined in terms of continued fraction expansion, which is closely related, as shown by Yoccoz, to the dynamical properties of the quadratic polynomial (see \cite{Y}). 
Classes of numbers defined by a condition analogous to \eqref{brjunowrtomega} when $d=2$, which is slightly stronger than Brjuno, were constructed in \cite{MMY}. 
Condition \eqref{brjuno-russmann} was introduced by Rssmann in KAM theory, making it possible to deal with a vector of frequencies. Brjuno-Rssmann conditions are already known to be central in the study of the linearization of vector fields (see e.g. \cite{GM}, \cite{P} and references therein). 

\bigskip
The classical Brjuno condition on the frequency was also considered by Young in \cite{Y93}, who constructed examples of non-reducible discrete cocycles in this case. For discrete cocycles with one frequency, Zhou and Wang recently obtained in \cite{ZW} a positive measure reducibility result for non-Brjuno frequencies for non-degenerate one-parameters families of cocycles.
Other results have been obtained on quasiperiodic cocycles regardless of any arithmetic condition on the frequency, worth mentioning although they are not reducibility results in our sense. 
In \cite{AFK}, it is shown that without any condition on the frequency, the Schrdinger cocycle \eqref{schr} can be conjugated to a rotation-valued cocycle for a positive measure set of energies; in \cite{ZW}, Zhou and Wang showed that in the case $d=2$, for any frequency, for a non-degenerate analytic one-parameter family which is close to a constant, the cocycle can be analytically diagonalized for a positive measure set of parameters.

\bigskip
\noindent
Our main result states:

\begin{thm}Let $\omega$ be a Brjuno vector, $A\in sl(2,\mathbb{R})$, $r>0$, $F\in C^\omega_r(\mathbb{T}^d,sl(2,\mathbb{R}))$. Suppose $\rho(A+F)$ is a $\frac{1}{2}$-Brjuno number with respect to $\omega$. There exists $\epsilon_0$ depending only on $\omega,\rho(A+F),r$ such that if $\mid F\mid_r\leq \epsilon_0$, then there exists $r'\in (0,r)$ such that $A+F$ is reducible in $C^\omega _{r'}$.
\end{thm}

\bigskip
\noindent
This result gives an extension of Eliasson's theorem using Rssmann's and Pschel's formulation of arithmetic conditions by means of approximation functions and their use in KAM methods. 
It holds for cocycles with arbitrarily many frequencies and gives a quite explicit link between reducibility and the arithmetics of $\omega$ and $\rho$, what Zhou-Wang's result does not since they consider larger dimensional systems with only two frequencies and formulate reducibility in terms of an abstract parameter. However, Zhou-Wang's article implies that it is impossible to find a lower bound for the function $(\omega,\rho)\mapsto \log \epsilon_0(\omega,\rho) + \int_1^\infty \frac{\log G(t)}{t^2}+\frac{\log g(t)}{t^{\frac{3}{2}}}dt$: Brjuno-Rssmann conditions are therefore not optimal in this problem as they might be in other dynamical problems.

\bigskip
\noindent In fact, our method gives this more explicit theorem:

\begin{thm}\label{main} Let $\kappa>0$ and let $G,g$ be positive increasing functions such that
\begin{itemize}
\item $G(1)\geq 1, g(1)\geq 1$,
\item

\begin{equation}\int_1^{+\infty}\frac{\log G(t)+\log g(t)}{t^2}dt<+\infty,\end{equation}

\item the map $t\mapsto \frac{g(t^2)}{G(t)}$ is bounded.
\end{itemize}

Suppose $\omega$ has $G$ as an approximation function with constant $\kappa$. Let $A\in sl(2,\mathbb{R})$, $r>0$, $F\in C^\omega_r(\mathbb{T}^d,sl(2,\mathbb{R}))$. Let $ n_0\in \mathbb{N}$. 
There exist $\epsilon_0$ depending only on $g,\kappa,G,n_0, r$ such that if 

\begin{enumerate}

\item 

$$\mid F\mid_r\leq \epsilon_0,$$ 

\item $\rho(A+F)$ has $g$ as an approximation function with respect to $\omega$ with constant $\kappa'>\kappa\sup_{t\geq n_0}
\frac{g(t^2)}{G(t)}$,

\end{enumerate}
%
then there exists $r'\in (0,r)$ such that $A+F$ is reducible in $C^\omega_{r'}(2\mathbb{T}^d,sl(2,\mathbb{R}))$.

\end{thm}

\bigskip
\noindent 
A discussion on the dependence of $\epsilon_0$ on $g$, $G$ and the other parameters is given in subsection \ref{Particular}.
%
%
%
%
%
%

\noindent As an application, we consider the case when $g$ and $G$ look like exponentials (section \ref{Particular}):

\begin{thm}\label{new} Let $\kappa>0,\kappa'>0$ and let $G(t)=e^{\frac{t}{(\log t)^\delta}},g(t)=e^{t^\alpha},\delta>1,\alpha<1$. Suppose $\omega$ has $G$ as an approximation function with constant $\kappa$. Let $A\in sl(2,\mathbb{R})$, $r>0$, $F\in C^\omega_r(\mathbb{T}^d,sl(2,\mathbb{R}))$. 
There exist $\epsilon_0$ depending only on $\alpha,\kappa,\delta,\kappa', r$ such that if 

\begin{enumerate}

\item 

$$\mid F\mid_r\leq \epsilon_0,$$ 

\item $\rho(A+F)$ has $g$ as an approximation function with respect to $\omega$ with constant $\kappa'$,

\end{enumerate}
%
then there exists $r'\in (0,r)$ such that $A+F$ is reducible in $C^\omega_{r'}(2\mathbb{T}^d,sl(2,\mathbb{R}))$.

\end{thm}

\noindent 
Our aim is to adapt the Pschel-Rssmann method (see \cite{R} and \cite{P}), which was used in the problem of linearization for vector fields, to the problem of reducing cocycles. It is a KAM-type method in which the speed of convergence is linear. 


\noindent 
First of all, we will build a setup in which a system $A+F$ with $A$ constant and $F$ small is conjugated to another system which is arbitrarily close to a constant, in an analytic class which, however, cannot be well controlled: this follows the technique used in \cite{El01} and is obtained by iterating (as in subsection \ref{iteration}) infinitely many steps (described in section \ref{step}) in which one conjugates a system $A_n+F_n$ to a system $A_{n+1}+F_{n+1}$ where $A_n,A_{n+1}$ are constant and $\mid F_{n+1}\mid_{r_n+1}\leq C \mid F_n\mid_{r_n}$, with $C<1$ being independent of $n$ and $r_n$ being a decreasing sequence controlling how analytic a function is. Thus, if $r_n$ tends to a non zero limit, we have analytic reducibility. 

\noindent 
At each step, in order to proceed, the constant part has to be non resonant, and if it is resonant, then we will have to remove the resonances, as explained in subsection \ref{elimination}.

\noindent 
Our setup makes sure that $r_n$ tends to a non zero limit whenever there is only a small enough number of steps at which one has to remove resonances in the constant part. The Brjuno-Rssmann condition on the frequency and on the rotation number of the cocycle is required exactly at this stage. 

\paragraph{Acknowledgments} The first author would like to thank the Centro Ennio de Giorgi for its hospitality and financial support during one year. Both authors would like to express their gratitude to H\aa kan Eliasson for useful discussions which allowed them to substantially improve the paper.

\section{The basic step}\label{step}

In this section, we will prove the iterative step, which consists in conjugating a system to another one with a smaller non constant part, whether the constant part be resonant or not.

\subsection{Definitions and notations}\label{defnot}

\noindent 
We will adopt the following conventions:

\bigskip

\noindent 
\Def Let $m=(m_1,\dots,m_d)\in\mathbb{Z}^d$. We denote by $\mid m\mid $ its modulus:\\
 $\mid m\mid= \sum_{j=1}^d \mid m_j\mid$.

\bigskip

\noindent 
\Def Let $F\in C^0(\mathbb{T}^d)$ and $r>0$; we say that $F\in C^\omega_r(\mathbb{T}^d)$ if there exists an analytic continuation of $F$ on a product of strips $\{(z_1,\dots, z_d)\in \mathbb{C}^d,\ \forall j\  \mid \operatorname{Im} z_j\mid <r \}$ and if the weighted norm

\begin{equation}\mid F\mid _r = \sum_{k\in \mathbb{Z}^d} \mid \mid \hat{F}(k)\mid \mid e^{2\pi \mid k\mid r}
\end{equation}

\noindent 
where $\mid \mid .\mid \mid $ is the relevant norm for $ \hat{F}(k)$ (for matrices, we use the operator norm), is finite.

\bigskip
\noindent Note that $C^\omega_r(\mathbb{T}^d)$ is a Banach space.

\bigskip

\noindent 
\textbf{Notation:} For $F\in C^\omega_r(\mathbb{T}^d)$, we denote its truncation by

$$F^N(\theta)=\sum_{\mid m\mid \leq N} \hat{F}(m)e^{2i\pi \langle m,\theta\rangle}.$$

\bigskip
\noindent 
\rem The weighted norms are particularly convenient since they satisfy, for any integer $N$,

\begin{equation}\mid F-F^N\mid_r=\sum_{k\in\mathbb{Z}^d,\mid k\mid > N}\mid \mid \hat{F}(k)\mid \mid e^{2\pi \mid k\mid r} =\mid F\mid_r-\mid F^N\mid_r  .
\end{equation}

\noindent Moreover, they are related (although not equivalent) to the usual sup norms since it is easy to see that

$$\sup_{\mid \operatorname{Im}\theta\mid  <r}\mid\mid F(\theta)\mid\mid 
\leq \mid F\mid_r.$$

\noindent For $r'<r$, we still have

$$\mid F\mid_{r'}\leq C(r-r')\sup_{\mid \operatorname{Im}\theta\mid  <r}\mid\mid F(\theta)\mid\mid $$

\noindent 
where $C(r-r')$ does not depend on $F$ but depends on $r-r'$. 

\bigskip
\noindent 
\Def If $A\in C^0(\mathbb{T}^d,sl(2,\mathbb{R}))$ and $X$ is the solution of $\frac{d}{dt}X(t,\theta)=A(\theta+t\omega)X(t,\theta);\ X(0,\theta)=Id$, the rotation number $\rho (A)$ is the quantity

$$\rho(A)= \lim_{t\rightarrow +\infty} \frac{1}{t} \operatorname{Arg} (X(t,\theta)\phi - \phi)$$

\noindent where $\phi \in \mathbb{R}^2\simeq \mathbb{C}$, $\theta\in \mathbb{T}^d$ and $\operatorname{Arg} $ stands for the variation of the complex argument ($\rho(A)$ is independent of $\theta, \phi$).

\bigskip
\noindent
Here and in what follows, we will fix $\omega\in \mathbb{R}^d$ rationally independent (i.e such that for all non zero $m\in \mathbb{Z}^d, \langle m,\omega\rangle\neq 0$): the vector $\omega$ will be the frequency of the cocycles we will consider.
We will always assume that $\kappa,\kappa'>0$ and that $G,g$ are two positive continuous and strictly increasing functions such that $G(1)\geq 1,g(1)\geq 1$.

\bigskip
\noindent 
\Def
$ NR(\kappa,G)=\{ \omega \in \mathbb{R}^d, \forall m\in\mathbb{Z}^d\setminus\{0\}, \mid \langle m,\omega\rangle\mid \geq \frac{\kappa}{G(\mid m\mid) }\}$.

\bigskip
\noindent
\rem There exists  a positive increasing and unbounded function $G\in C^0(\mathbb{R}^{*+})$ with $G(1)\geq 1$ and $\kappa>0$ such that $\omega \in \mathrm{NR}(\kappa,G)$. 

\noindent
Indeed, one can take $\kappa=\min_i \mid \omega_i\mid $ and $G(N)=\max_{\mid m\mid \leq N}
\frac{\kappa}{\mid \langle m,\omega \rangle\mid}$.

\noindent 
As noticed by H. Rssmann (\cite{R}), a condition $\mathrm{NR}(\kappa,G)$ with $G$ such that 

\begin{equation}\label{condomega}\int_1^\infty\frac{\log G(t)}{t^2}dt<\infty
\end{equation}

\noindent is fulfilled by all Bruno vectors (see \cite{B}), i.e vectors satisfying:

\begin{equation}\label{br}\sum_{k\geq 1} \frac{\mid \log \alpha_{2^k-1}\mid }{2^k}<+\infty
\end{equation}

\noindent 
where 

\begin{equation}\alpha_k=\min_{l\leq k}\min_{j=1,\dots,d}\min_{\mid m\mid =l+1}\mid \langle m,\omega\rangle 
-\omega_j\mid .
\end{equation}

\noindent 
In \cite{GM} it is shown that condition \eqref{br} is equivalent to

\begin{equation}\label{tau}\sum_{k\geq 1} \frac{\mid \log \alpha_k\mid}{k(k+1)}<+\infty
\end{equation}

\noindent 
and condition \eqref{condomega} implies condition \eqref{tau}, so that \eqref{condomega}, \eqref{br} and \eqref{tau} are equivalent. 
Moreover, in dimension $d=2$, these conditions are equivalent to the usual Brjuno condition on $\frac{\omega_2}{\omega_1}$ (see \cite{Y}). This suggests the following definition:

\bigskip
\noindent \Def The vector $\omega$ is a Brjuno vector if $\omega\in\mathrm{NR}(\kappa,G)$ with $G$ satisfying \eqref{condomega}.

\bigskip
\noindent We now have to introduce another type of arithmetic condition, related to the well-known "second Melnikov condition".

\bigskip
\noindent 
\Def Let $N\in \mathbb{N}\setminus\{0\}$; we set

\begin{equation}\mathrm{NR}_\omega^N(\kappa', g)=\left\{\alpha\in \mathbb{C},
\forall m\in \mathbb{Z}^d\setminus\{0\}, 0<\mid m\mid \leq N\Rightarrow \mid \alpha-i\pi \langle m,\omega\rangle \mid \geq \frac{\kappa'}{g(\mid m\mid )}\right\}
\end{equation}

\noindent 
and $\mathrm{NR}_\omega(\kappa', g)=\cap_{N\in \mathbb{N}} \mathrm{NR}_\omega ^N (\kappa',g)$. 

\bigskip
\noindent \rem If $g(t)=t^\tau$ for some $\tau>1$, this is a diophantine condition. 

\bigskip
\noindent \Def Let $\nu>0$. The number $\alpha$ is a $\nu$-Brjuno number with respect to $\omega$ if $\alpha\in \mathrm{NR}_\omega (\kappa',g)$ with $g$ satisfying

\begin{equation}\label{gnu}\int_1^\infty \frac{\log g(t)}{t^{(1+\nu)}}<+\infty.
\end{equation}

\noindent This extended Brjuno condition was first considered in \cite{MMY}.

\subsection{Elimination of resonances}\label{elimination}

\noindent 
Now we shall prove the uniqueness of a resonance, i.e the situation of the spectrum being close to a number of the form $\langle m,\omega\rangle,m\in \mathbb{Z}^d$, when it exists.

\begin{lem}\label{m} Let $\alpha\in \mathbb{C}$. Let $N\in \mathbb{N}\setminus\{0\} $.
There exists $m\in \mathbb{Z}^d$ such that $\mid m\mid \leq N$ and $\alpha-i\pi\langle m,\omega\rangle \in 
\mathrm{NR}_\omega^N(\frac{\kappa}{4G(N)},g)$; 
if $m$ is non zero, then 

$$\mid \alpha-i\pi \langle m,\omega\rangle \mid < \frac{\kappa}{4G(N)g(\mid m\mid )}$$

\noindent and $\alpha-i\pi \langle m,\omega\rangle\in \mathrm{NR}_\omega^N(\frac{\kappa}{G(N)},g)$. 

\end{lem}

\noindent 
\dem 
Suppose $\alpha$ is not in $\mathrm{NR}_\omega^N(\frac{\kappa}{4G(N)},g)$, i.e there exists $m\in \mathbb{Z}^d, 0<\mid m\mid \leq N$, such that 

$$\mid \alpha-i\pi \langle m,\omega\rangle \mid < \frac{\kappa}{4G(N)g(\mid m\mid)}.$$ 

\noindent 
Then for all $m'\in \mathbb{Z}^d$ with $0<\mid m'\mid \leq N$, 

\begin{equation}\mid  \alpha-i\pi \langle m+m',\omega\rangle \mid \geq \mid \pi \langle m',\omega\rangle \mid -\mid \alpha-i\pi \langle m,\omega\rangle \mid
\geq \frac{\kappa}{G(\mid m'\mid )}-\frac{\kappa}{4G(N)g(\mid m\mid)}
\end{equation}

\noindent 
so

\begin{equation}\mid  \alpha-i\pi \langle m+m',\omega\rangle \mid \geq \frac{\kappa}{G(N)g(\mid m'\mid)}
\end{equation}

\noindent 
and so $\alpha-i\pi \langle m,\omega\rangle\in \mathrm{NR}_\omega^N(\frac{\kappa}{G(N)},g)$. $\Box$

\bigskip

\noindent 
The following proposition explains how to eliminate resonances in the spectrum of a trace zero matrix.

\begin{prop}\label{phi} Let $A\in sl(2,\mathbb{R})$ with eigenvalues $\pm \alpha$. Let $N\in \mathbb{N}$.
Suppose that $\alpha$ is not in $\mathrm{NR}_\omega^N(\frac{\kappa}{4G(N)},g)$.
There exists $\Phi \in \cap_{r'\geq 0} C^\omega_{r'}(2\mathbb{T}^d,SL(2,\mathbb{R}))$ and a numerical constant $C'$ such that

\begin{equation}\label{Phi}
\forall r'\geq 0, \ \mid \Phi\mid_{r'} \leq C'e^{\pi Nr'}; \  \mid \Phi^{-1}\mid_{r'} \leq C'e^{\pi Nr'}
\end{equation}

\noindent 
and if $\tilde{A} $ with eigenvalues $\pm \tilde{\alpha}$ is such that

\begin{equation}\label{tildeA} \partial_\omega \Phi=A\Phi - \Phi \tilde{A}
\end{equation}

\noindent 
then $\tilde{\alpha} \in \mathrm{NR}_\omega^{N} (\frac{\kappa}{G(N)},g)$. 


\noindent 
Moreover,

$$\mid \tilde{\alpha} \mid < \frac{\kappa}{4G(N) }.$$

\end{prop}

\dem Lemma \ref{m} gives a number $m, 0<\mid m\mid \leq N$, such that letting 

$$\tilde{\alpha}=\alpha-i\pi \langle m,\omega\rangle$$ 

\noindent 
then $\tilde{\alpha}\in \mathrm{NR}_\omega^N(\frac{\kappa}{G(N)},g)$.
Let $P$ be such that $P^{-1}AP$ is diagonal and $\mid \mid P\mid \mid =1$. We define 

$$\Phi(\theta)=P^{-1}\left(\begin{array}{cc}
e^{i\pi \langle m,\theta\rangle} & 0\\
0 & e^{-i\pi \langle m,\theta\rangle}\\
\end{array}\right)P.$$

\noindent 
Relation \eqref{tildeA} gives

$$\tilde{A}=P^{-1}\left(\begin{array}{cc}
\tilde{\alpha} & 0\\
0 & -\tilde{\alpha}\\
\end{array}\right)P.
$$

\noindent 
To obtain the estimate \eqref{Phi}, we use an estimate shown for instance in \cite{El01}, Lemma A:

$$\mid \mid P^{-1}\mid \mid \leq \max \left(1,\left(\frac{C.\mid\mid A\mid\mid}{2\mid\alpha\mid}\right)^6\right)
$$

\noindent where $C$ is a numerical constant,
and since $A$ is diagonalizable whenever its eigenvalues are non zero, 

$$\mid \mid P^{-1}\mid \mid \leq \max \left(1,\left(\frac{C.\mid\alpha\mid}{2\mid\alpha\mid}\right)^6\right)
\leq C'
$$

\noindent 
where $C'$ is a numerical constant, 
which  gives \eqref{Phi}. $\Box$

\subsection{Solution of the linearized homological equation}\label{KAM}

\noindent 
Our aim is to solve an equation of the form

$$\partial_\omega Z=(A+F)Z-Z(A'+F')$$

\noindent 
where $A$ and $F$ are known, $A\in sl(2,\mathbb{R})$ and $F$ is analytic with values in $sl(2,\mathbb{R})$. If $A$ is non-resonant, we first solve

$$\partial_\omega \tilde{X}=[A,\tilde{X}]+aF^N-a\hat{F}(0)$$

\noindent 
where $F^N$ is some truncation of $F$ and $a$ is close enough to 1; then we define $A'=A+a\hat{F}(0)$ and $F'$ by

$$\partial_\omega e^{\tilde{X}}=(A+F)e^{\tilde{X}}-e^{\tilde{X}}(A'+F')$$

\noindent 
and then we estimate $F'$ to get $\mid F'\mid_{r'}\leq \sqrt{1-a}\mid F\mid_r$. If $A$ is resonant, we conjugate $A+F$ to a  system $\tilde{A}+\tilde{F}$ where $\tilde{A}$ is non-resonant and we proceed in the same way as in the non-resonant case.  So, from now on, to simplify the notations, we will assume that $A,\tilde{A},A'\in sl(2,\mathbb{R})$, that $F,\tilde{F},\tilde{X},F'$ have values in $sl(2,\mathbb{R})$ and $Z,\Phi$ have values in $SL(2,\mathbb{R})$.

\begin{prop}\label{solhomol}
Let $N\in \mathbb{N}$ and $r,r'>0$.

\noindent 
Let $\tilde{A}$ with eigenvalues $\pm \tilde{\alpha}\in \mathrm{NR}_\omega^{N}(\frac{\kappa}{4G(N)},g)$. Let $\tilde{F}
\in C^\omega _r(\mathbb{T}^d)$. Then equation

\begin{equation}\label{homol}
\forall \theta\in \mathbb{T}^d, \ \partial_\omega \tilde{X}(\theta)= [ \tilde{A}, \tilde{X}(\theta)] +\tilde{F}^{N}(\theta)-\hat{\tilde{F}}(0);\ \hat{\tilde{X}}(0)=0
\end{equation}

\noindent 
has a unique solution $\tilde{X}\in C^\omega_{r'}(\mathbb{T}^d)$ such that
 
 \begin{equation}
  \mid  \tilde{X}\mid _{r'}\leq \frac{4}{\kappa}G(N)g(N)\mid \tilde{F}^{N}\mid _{r'}
 \end{equation} \end{prop}
 
 
 \dem In Fourier series, equation \eqref{homol} can be written:
 
 \begin{equation}\label{homolfou}\begin{split} \forall m\in \mathbb{Z}^d, &0<\mid m\mid \leq N \Rightarrow
 2i\pi \langle m,\omega\rangle \hat{\tilde{X}}(m)= [\tilde{A}, \hat{\tilde{X}}(m)] +\hat{\tilde{F}}(m);\\
& \mid m\mid \in \{0\}\cup [N+1,+\infty[ \Rightarrow 2i\pi \langle m,\omega\rangle \hat{\tilde{X}}(m)= [\tilde{A}, \hat{\tilde{X}}(m)] .
\end{split} \end{equation}
 
\noindent 
 So for $ \mid m\mid \in \{0\}\cup [N+1,+\infty[$, $ \hat{\tilde{X}}(m)=0$ is a solution (not necessarily unique). 
 
\noindent 
 For $0<\mid m\mid \leq N$, 
 the solution is formally written as
 
 \begin{equation}\hat{\tilde{X}}(m)=\mathcal{L}_m^{-1}\hat{\tilde{F}}(m)
 \end{equation}
 
\noindent 
 where $\mathcal{L}_m$ is the operator 
 
 $$\mathcal{L}_m:sl(2,\mathbb{R})\rightarrow sl(2,\mathbb{R}),\ M\mapsto 2i\pi \langle m,\omega\rangle M - [\tilde{A},M].$$

\noindent 
 Its spectrum is 
 $\{2i\pi \langle m,\omega\rangle-2\tilde{\alpha},2i\pi \langle m,\omega\rangle+2\tilde{\alpha},2i\pi \langle m,\omega\rangle \}$.

\noindent 
 Since $\omega\in \mathrm{NR}(\kappa,G)$ and $\tilde{\alpha}\in \mathrm{NR}_\omega^{N}(\frac{\kappa}{4G(N)},g)$, $\mathcal{L}_m$ is invertible and we have for all $m\in\mathbb{Z}^d$ such that $\mid m\mid\in (0, N]$,
 
 $$\mid \mid \mathcal{L}_m^{-1}\mid \mid \leq \max\{\frac{G(\mid m\mid )}{\kappa},
 \frac{4G(N)g(\mid m\mid )}{\kappa}\}= \frac{4G(N)g(\mid m\mid )}{\kappa}$$ 
 
\noindent 
 therefore for all $m\in \mathbb{Z}^d$ such that $0<\mid m\mid \leq N$,
 
 \begin{equation}
 \mid \mid \hat{\tilde{X}}(m)\mid \mid \leq   4G(N) \frac{g(\mid m\mid) }{\kappa}\mid \mid \hat{\tilde{F}}(m)\mid \mid 
 \end{equation}
 
\noindent 
  therefore
 
 \begin{equation}\begin{split}
 \mid  \tilde{X}\mid _{r'}&\leq 4G(N)\sum_{m\in \mathbb{Z}^d\setminus\{0\},\mid m\mid\leq N} \frac{g(\mid m\mid)}{\kappa} \mid \mid \hat{\tilde{F}}(m)\mid \mid 
e^{2\pi\mid m\mid r'}\\
& \leq 4\frac{G(N)g(N)}{\kappa} \mid \tilde{F}^{N}\mid _{r'}.\ \Box
 \end{split}\end{equation}

 \subsection{Solution of the full homological equation without resonances}
 
\noindent 
 This section explains the basic step in case the constant part is non resonant, i.e when its eigenvalues are far from all $\langle m,\omega\rangle$, $m\in\mathbb{Z}^d\setminus\{0\}$.
 
 \begin{prop}\label{iternonres}
 Let $0<r'\leq r$, $a'\in (0,1]$, $N\in \mathbb{N}$, $\tilde{F} \in C^\omega_r(\mathbb{T}^d)$, $\tilde{A}\in sl(2,\mathbb{R})$. If $\sigma(\tilde{A})=\{\pm \alpha\}, \alpha\in \mathrm{NR}_\omega^{N}(\frac{\kappa}{4G(N)},g)$,  then there exists $\tilde{X},F'\in C^\omega_{r'}(\mathbb{T}^d)$and $A'\in sl(2,\mathbb{R})$
 such that

 \begin{equation}\label{a'close}\mid A'-\tilde{A}\mid \leq \mid \mid \hat{\tilde{F}}(0)\mid \mid
 \end{equation}

 \begin{equation}\label{conj}\partial_\omega e^{\tilde{X}}=(\tilde{A}+\tilde{F})e^{\tilde{X}}-e^{\tilde{X}}(A'+F')
 \end{equation}


 \begin{equation}\label{Xtilde}
 \mid \tilde{X} \mid _{r'}\leq 4 a' \frac{G(N)g(N)}{\kappa}\mid \tilde{F}^{N}\mid _{r'}
 \end{equation}

\noindent 
and

\begin{equation}\label{F'0}\begin{split}
\mid F' \mid_{r'} &\leq e^{\mid \tilde{X} \mid_{r'}}(1-a')\mid  \tilde{F} \mid _{r'}+e^{\mid \tilde{X} \mid_{r'}}a'\mid \tilde{F}-\tilde{F}^{N}\mid_{r'}\\
&+e^{\mid \tilde{X} \mid_{r'}} \mid  \tilde{F} \mid _{r'}\mid \tilde{X}\mid_{r'}(e^{\mid \tilde{X} \mid_{r'}}+a'+a'e^{\mid \tilde{X} \mid_{r'}}).
\end{split}\end{equation}

\end{prop}

\dem Let $\tilde{X}$ be a solution of

\begin{equation}\label{homola}
\forall \theta\in \mathbb{T}^d, \ \partial_\omega \tilde{X}(\theta)= [ \tilde{A}, \tilde{X}(\theta)] +a'\tilde{F}^{N}(\theta)-a'\hat{\tilde{F}}(0);\ \hat{\tilde{X}}(0)=0
\end{equation}

\noindent as given by Proposition \ref{solhomol} (so it satisfies \eqref{Xtilde}). Let $A'=\tilde{A}+a'\hat{\tilde{F}}(0)$ so that \eqref{a'close} holds, and let $F'$ be defined by 

$$\partial_\omega e^{\tilde{X}}=(\tilde{A}+\tilde{F})e^{\tilde{X}}-e^{\tilde{X}}(A'+F').$$

\noindent 
We have

\begin{equation}\begin{split}F'=&e^{-\tilde{X}} (\tilde{F}-a'\tilde{F}^{N})+e^{-\tilde{X}}\tilde{F}(e^{\tilde{X}}-Id)+a'(e^{-\tilde{X}}-Id)\hat{\tilde{F}}(0)\\
&-e^{-\tilde{X}}\sum_{k\geq 2}\frac{1}{k!}\sum_{l=0}^{k-1}\tilde{X}^l(a'\tilde{F}^{N}-a'\hat{\tilde{F}}(0))\tilde{X}^{k-1-l}.
\end{split}\end{equation}

\noindent 
Since

\begin{equation} \mid \tilde{F}-a'\tilde{F}^{N} \mid_{r'}\leq 
a'\mid \tilde{F}-\tilde{F}^{N}\mid_{r'}+(1-a')\mid \tilde{F} \mid _{r'}.
\end{equation}

\noindent 
one easily obtains \eqref{F'0}. $\Box$
%
%

\bigskip
\rem Denote $\epsilon= \mid \tilde{F}\mid_r$. 
Suppose

\begin{equation}\label{petitesse3}  2G(N)g(N)\epsilon\leq \frac{  \kappa(1-a')}{2}.
\end{equation}

\noindent 
Then \eqref{Xtilde} implies $\nonumber \mid \tilde{X}\mid_{r'}\leq a'(1-a')
$
thus $e^{\mid \tilde{X}\mid_{r'}}\leq 2$.
By \eqref{F'0}, if one assumes moreover that
%


\begin{equation}e^{-2\pi N(r-r')}\leq 1-a'\end{equation}

\noindent 
with $r'>0$,
then 

$$\mid \tilde{F}-\tilde{F}^{N}\mid_{r'}\leq (1-a')\mid \tilde{F}-\tilde{F}^{N}\mid_{r}$$

\noindent 
therefore

\begin{equation}\mid F'\mid _{r'}\leq 
2(1-a')\epsilon+ 2a' (1-a')\epsilon+2\epsilon a'(1-a')(2+3a').
\end{equation}

\noindent 
Thus, if $a'$ is close enough to 1 (i.e larger than $1-\frac{1}{14^2}$), 

\begin{equation}\mid F'\mid _{r'}\leq (1-a')^{\frac{1}{2}}\epsilon .\end{equation}

 \subsection{Solution of the full homological equation with resonances}
 
\noindent 
 This section presents the basic step when there are resonances in the constant part, i.e when its eigenvalues are too close to some $\langle m,\omega\rangle, m\in \mathbb{Z}^d\setminus\{0\}$.

 \begin{prop}\label{iter2} 
Let $a\in (0,1),c_0>0$, $C'$ as in Proposition \ref{phi}, $N\in \mathbb{N}$, $r>\frac{2\log (g(N)G(N))}{\pi N} $, ${F} \in C^\omega_r(\mathbb{T}^d )$, $\mid F\mid_r=\epsilon$,
 $A\in sl(2,\mathbb{R})$. Suppose the eigenvalues $\pm \alpha$ of 
 $A$ are not in $\mathrm{NR}_\omega^N(\frac{\kappa}{4G(N)},g)$. If
 
 \begin{equation}\label{petitesse}
2 G(N)^2g(N)^2\epsilon\leq \frac{(1-a)^2}{2 }
 \kappa^2
 \end{equation}
\noindent 
and 

\begin{equation}\label{petitesse2}eC'(G\cdot g)(N+1)^{-c_0}   \leq 1-a
\end{equation}


\noindent 
then letting 
$r'=\frac{r}{2}-c_0\frac{\log (G\cdot g)(N+1)}{4\pi N}$, there exists ${F}'\in C^\omega_{r'}(\mathbb{T}^d )$,
$A'\in sl(2,\mathbb{R})$ and $Z\in  C^\omega_{r'}(2\mathbb{T}^d )$ such that


 \begin{equation}\label{Z}\partial_\omega Z=(A+F)Z-Z (A'+F')
 \end{equation}
 
 \noindent and

 \begin{equation}\label{F'}\mid F'\mid _{r'}\leq  (1-a) \epsilon.
 \end{equation}

 \end{prop}

\dem One first applies Proposition \ref{phi} on $A$. Let $\Phi,\tilde{A}$ be as in Proposition \ref{phi} so that   
$\sigma(\tilde{A})=\pm \tilde{\alpha}$, 
$\tilde{\alpha}\in \mathrm{NR}_\omega^N(\frac{\kappa}{G(N)},g)$ and $\mid\tilde{\alpha}\mid \leq \frac{\kappa}{4G(N)}$; let $\tilde{F}= \Phi F \Phi^{-1}$. Notice that, by construction of $\Phi$, the map $\tilde{F}$ remains continuous on 
$\mathbb{T}^d$. 
Apply Proposition \ref{iternonres} with $a'=1$ and with $r'=\frac{r}{2}-c_0\frac{\log (G\cdot g)(N+1)}{4\pi N}$ to get $\tilde{X}\in C^\omega_{r'}(\mathbb{T}^d ), A',F'\in C^\omega_{r'}(\mathbb{T}^d )$ such that \eqref{conj} and \eqref{F'0} hold as well as

 \begin{equation}
 \mid \tilde{X} \mid _{r'}\leq  4\frac{G(N)g(N)}{\kappa}\mid \tilde{F}^{N}\mid _{r'}
 \end{equation}

\noindent 
and let $Z=\Phi e^{\tilde{X}}\in C^\omega_{r'}(2\mathbb{T}^d )$ so that $Z$ satisfies \eqref{Z}. 
Condition \eqref{petitesse} 
implies that 

$$(G\cdot g)(N)\mid \tilde{X}\mid _{r'}\leq   \frac{(1-a)^2\mid \tilde{F}^{N}\mid_{r'}}{\epsilon  }$$ 

\noindent 
so \eqref{F'0} with $a'=1$ gives

\begin{equation}\mid F'\mid_{r'}\leq eC'\mid F-F^N\mid_{r}e^{-2\pi N(r-2r')}   +eC'\mid F^N\mid_{r'}e^{2\pi Nr'}  \frac{(1-a)^2}{ (G\cdot g)(N)}(2e+1)\end{equation}

\noindent and by the choice of $r'$,

\begin{equation}\mid F'\mid_{r'}
\leq eC'\mid F\mid_{r}(G\cdot g)(N+1)^{-c_0}.   \end{equation}

\noindent 
This implies, by assumption \eqref{petitesse2}, that

\begin{equation}\begin{split}\mid F' \mid_{r'} \leq (1-a)\mid F\mid_r.\ \Box\end{split}\end{equation}

\section{Iteration, reducibility and arithmetical conditions}
\subsection{Iteration}\label{iteration}

\noindent 
In this section we will first introduce the Brjuno-Rssmann condition and then show how one can use it to control the convergence of the KAM iteration scheme.

\begin{ass}\label{brunoGg} The functions $g$ and $G$ satisfy

\begin{equation}\label{condomega'}\int_1^\infty \frac{\log [g(t)G(t)]}{t^2}dt<\infty.
\end{equation}

\end{ass}

\noindent 
In order to iterate the basic step, we will now fix the parameters as follows: let $C'$ be as in Proposition \ref{phi}. Furthermore, let $r_0>0$, $n_0\in \mathbb{N}$ and choose

$$c_0=\frac{r_0}{4^{n_0+3}(\sup_{t\in [1,n_0]} \frac{\log (G\cdot g)(t+1)}{t}+1)}$$

\noindent Let $a\in [1-\bar{a},1)$ where $\bar{a}=\min (\frac{1}{14^2}, \frac{1}{(G\cdot g)(2)^2})$.
Let $\epsilon_0>0$ be small enough to assure that

\begin{equation}\label{condepsilon}\int_{(G\cdot g)^{-1}\left(\frac{\kappa }{2(1-a)^{\frac{n_0-5}{4}}\epsilon_0^{\frac{1}{2}}}\right)}^\infty \frac{\log (G\cdot g)(t)}{t^2}dt\leq \frac{r_0}{4^{n_0+2}}
\end{equation}

\noindent 
and  

\begin{equation} \label{condepsilon2}e C' \epsilon_0^{\frac{c_0}{4}}\leq  (1-a)^2\kappa^2.
\end{equation}

\noindent For all $n\in \mathbb{N}$, let $\epsilon_n=(1-a)^{\frac{n}{2}}\epsilon_0$ and let $N_n$ be the biggest integer such that
%
$$
(G\cdot g)(N_n) ^2
\leq \frac{(1-a)^2}{4 \epsilon_n}\kappa^2$$

\noindent ($N_n$ exists since $\epsilon_n\leq \frac{(1-a)^2\kappa^2}{4e (G\cdot g)(1)^2}$). 

\bigskip
\noindent The above choices of the sequences $\epsilon_n$ and $N_n$ are made in such a way that the following holds:
%
%
%
%
%
%

\begin{equation}
\label{Nn0}\int_{N_{n_0}}^\infty \frac{\log (G\cdot g)(t)}{t^2}dt\leq \frac{r_0}{4^{n_0+2}}.
\end{equation}




\noindent
\rem: The number $\epsilon_0$ will then only depend on $a,\kappa,g,G,n_0$ and $r_0$ (the larger $n_0$ and the smaller $r_0$ are, the smaller $\epsilon_0$ will be). 

\bigskip
\noindent 
To simplify the notations, from now on the functions $A_n$ are understood to be in $sl(2,\mathbb{R})$, while the $F_n$ have their values in $sl(2,\mathbb{R})$ and $Z_n',Z_n$ have their values in $SL(2,\mathbb{R})$.

\begin{prop}\label{red0}Let $A\in sl(2,\mathbb{R}) $ and $F\in C^\omega_{r_0}(\mathbb{T}^d )$. 
If $\mid F\mid_{r_0}\leq \epsilon_0$,
then there exist sequences $(r_n)_{n\in \mathbb{N}}, r_n>0$,
$Z_n,
\in C^\omega_{r_n}(2\mathbb{T}^d )$, 
$A_n $ with spectrum $\pm \alpha_n$, 
$F_n\in C^\omega_{r_n}(2\mathbb{T}^d )$,  
and $m_n\in\mathbb{Z}^d$,
such that 


\begin{enumerate}
\item \label{1} if all $m_n$ are zero when $n\geq n_0$, then $r_n$ has a positive limit;
\item \label{2} if $m_{n}\neq 0$ then $\mid \alpha_{n}-\pi \langle m_{n},\omega\rangle\mid \leq\frac{\kappa}{4G(N_n)}$; 
\item \label{3} $m_n$ has modulus  less than $N_{n}$,

\item \label{4} $\mid F_n\mid_{r_n}\leq \epsilon_n$;
\item \label{5} $\partial_\omega Z_n=(A+F)Z_n-Z_n(A_n+F_n);$

\item \label{7} $\mid \alpha_{n-1}-i\pi \langle m_{n-1},\omega\rangle -\alpha_n\mid \leq \sqrt{\epsilon_{n-1}}$.
\end{enumerate}

\end{prop}

\rem Proposition \ref{red0} implies that $A+F$ is reducible in $C^\omega_{r'}$ for some $r'>0$ if all $m_n$ are zero for $n\geq n_0$.


\bigskip
\noindent 
\dem This proposition is shown by recurrence. Suppose these sequences are defined up to some $n\in \mathbb{N}$ and suppose that for all $n'\leq \min(n-1,n_0)$, $r_{n'+1}\geq \frac{r_{n'}}{4}$. We must distinguish two cases according to the possibility that the spectrum of $A_n$ is resonant or not.

\bigskip
\textbf{First case:} 
$\alpha_n\in \mathrm{NR}_\omega^{N_n}(\frac{\kappa}{4G(N_n)},g)$. 
Let $r_{n+1}=r_n-c_0\frac{\mid \log (1-a)\mid}{2\pi N_n}$, so that $r_{n+1}\geq \frac{r_n}{2}$ if $n\leq n_0$.
One can apply Proposition \ref{iternonres} with $r=r_n$, $r'=r_{n+1}$, $N=N_n$,
$\tilde{F}=F_n$ and
$\tilde{A}=A_n$
and obtain $Z'_n=e^{\tilde{X}_n}\in C^\omega_{r_{n+1}}
(\mathbb{T}^d )$, $F_{n+1}\in C^\omega_{r_{n+1}}
(\mathbb{T}^d ), A_{n+1} $ such that 

\begin{equation}\partial_\omega Z'_n =(A_n+F_n)Z'_n -Z'_n (A_{n+1}+F_{n+1})
\end{equation}

\noindent and


 \begin{equation}\mid F_{n+1}\mid _{r_{n+1}}\leq (1-a)^{\frac{1}{2}}\epsilon_{n}=\epsilon_{n+1}.
 \end{equation}

\noindent 
One then takes $Z_{n+1}=Z_nZ'_n$.

\bigskip
\noindent 
\textbf{Second case:} $\alpha_n\notin \mathrm{NR}_\omega^{N_n}(\frac{\kappa}{4G(N_n)},g)$. 

\noindent 
Assumption \eqref{petitesse} is satisfied by definition of $N_n$; 
assumption \eqref{petitesse2} is also satisfied since, by maximality of $N_n$, 

\begin{equation}(G\cdot g)(N_n+1)^{-c_0}\leq \left(\frac{2\epsilon_n}{(1-a)^2\kappa^2}\right)^{\frac{c_0}{2}}
\leq  \left(\frac{2\epsilon_0}{(1-a)^2\kappa^2}\right)^{\frac{c_0}{2}}
\end{equation}

\noindent 
which, together with \eqref{condepsilon2}, implies that

\begin{equation}(G\cdot g)(N_n+1)^{-c_0}\leq \frac{1-a}{e C'}.
\end{equation}

\noindent 
Therefore, one can apply Proposition \ref{iter2} with $r=r_n$,
$r'=r_{n+1}=\frac{r_n}{2}-c_0\frac{\log (G\cdot g)(N_n+1)}{\pi N_n}$, so that $r_{n+1}\geq \frac{r_n}{4}$ if $n\leq n_0$,
 and $
N=N_n$. It follows that there exists $A_{n+1}\in sl(2,\mathbb{R})$,
$F_{n+1}\in C^\omega_{r_{n+1}}(\mathbb{T}^d )$ and $Z'_n\in C^\omega_{r_{n+1}}(2\mathbb{T}^d )$ such that 

$$\partial_\omega Z'_n=(A_n+F_n)Z'_n-Z'_n(A_{n+1}+F_{n+1})$$

\noindent and

$$\mid F_{n+1}\mid _{r_{n+1}}\leq (1-a) \mid F_{n}\mid _{r_{n}}\leq \epsilon_{n+1}
$$

\bigskip
\noindent 
One then takes $Z_{n+1}=Z_n
Z'_n$. 

\bigskip
\noindent 
To complete the proof, we now need to show that $(r_n)_n$ has a positive limit if all $m_n$ are zero for $n\geq N_0$. We have  

\begin{equation}\begin{split}\lim_n r_n=r_{n_0}-\sum_{k=n_0}^\infty (r_k-r_{k+1})
\geq \frac{r_0}{4^{n_0}}
- \sum_{k\geq  n_0} \frac{\mid \log (1-a)\mid}{2\pi N_k}.\\
\end{split}
\end{equation}

\noindent 
Now, for all $n$,

$$N_n=E\left((G\cdot g)^{-1}\left( \frac{(1-a)\kappa}{2\sqrt{ \epsilon_n}} \right) \right)$$


\noindent 
thus

\begin{equation}\begin{split}\lim_n r_n\geq   \frac{r_0}{4^{n_0}}
-\frac{\mid \log (1-a)\mid}{2\pi }\int_{n_0}^\infty \left[(G\cdot g)^{-1}\left({\frac{\kappa}{2 (1-a)^{\frac{n}{4}-1}}} \right)\right]^{-1}dn.
\end{split}\end{equation}

\noindent 
Through the change of variables $X= \frac{\kappa}{2 (1-a)^{\frac{n}{4}-1}} $,

\begin{equation}\begin{split}\lim_n r_n\geq   \frac{r_0}{4^{n_0}}
-\int_{(G\cdot g)(N_{n_0})}^\infty \frac{1}{\pi (G\cdot g)^{-1}(X) X}dX.
\end{split}\end{equation}

\noindent 
Letting now $Y=(G\cdot g)^{-1}(X)$, the integral becomes

\begin{equation}\begin{split}\int_{(G\cdot g)(N_{n_0})}^\infty \frac{1}{\pi (G\cdot g)^{-1}(X) X}dX&=\int_{N_{n_0}}^\infty \frac{1}{\pi Y (G\cdot g)(Y)}d(G\cdot g)(Y)\\
&=\frac{-\log (G\cdot g)(N_{n_0}) }{N_{n_0}}+\int_{N_{n_0}}^\infty \frac{\log (G\cdot g)(Y)}{Y^2 }dY.
\end{split}\end{equation}

\noindent 
Therefore 

\begin{equation}\lim_n r_n\geq   \frac{r_0}{4^{n_0}}+\frac{\log (G\cdot g)(N_{n_0}) }{\pi  N_{n_0}}-\frac{1}{\pi}\int_{N_{n_0}}^\infty \frac{\log (G\cdot g)(Y)}{Y^2 }dY
\end{equation}

\noindent 
so, by \eqref{Nn0},

\begin{equation}\lim_n r_n\geq  \frac{ r_0}{4^{n_0+1}}
\end{equation}

\noindent 
which is positive. $\Box$

\subsection{A link between the Brjuno sum and the allowed perturbation}

\noindent 
It is easily seen that the condition on $\epsilon_0$ can also be expressed more conveniently as the following sufficient condition:

\begin{equation}\label{brunosum}\epsilon_0\leq \exp \left( -\frac{r_0}{4^{n_0}}
-\mid \log \frac{\kappa}{2(1-a)^{n_0}}\mid -2\int_1^\infty \frac{\log (g\cdot G)(t)}{t^2}dt \right)
\end{equation}

\noindent Indeed, we have the following bound:

\begin{equation}\begin{split}\mid \frac{1}{2} \log \epsilon_0 +\int_1^\infty \frac{\log (g\cdot G)(t)}{t^2}dt)\mid
&=\mid \int_1^\infty \frac{\log (\sqrt{\epsilon_0} g\cdot G)(t)}{t^2}dt)\mid\\
&\leq \frac{r_0}{4^{n_0}}+\int_1^{(g\cdot G)^{-1}(\frac{\kappa}{2(1-a)^{n_0}\sqrt{\epsilon_0}})} \frac{\mid \log \sqrt{\epsilon_0} (g\cdot G)(t)\mid }{t^2}dt
\end{split}\end{equation}

\noindent and the conclusion follows easily from the upper bound on $t$.

\subsection{Reducibility theorem}

\noindent We will now need one more assumption on the approximation functions $G$ and $g$.
%
%
%

\begin{ass}\label{rapportgG} The map $t\mapsto \frac{g(t^2)}{G(t)}$ is bounded. \end{ass}

\noindent 
Now we can prove the main result:

\begin{thm}\label{redu} Let $A\in sl(2,\mathbb{R})$, $r>0$, $F\in C^\omega_r(\mathbb{T}^d )$. 
Let $ n_0\in \mathbb{N}$. 
Assume $\rho(A+F)\in \mathrm{NR}_\omega(\kappa',g)$ with $\kappa'>\kappa\sup_{t\geq n_0}
\frac{g(t^2)}{G(t)}$.

\noindent 
Under assumptions \ref{brunoGg} and \ref{rapportgG} on the approximation functions $g$ and $G$, there exist $\epsilon_0>0$ depending only on $g,\kappa,G,n_0, r$ such that if 

$$\mid F\mid_r\leq \epsilon_0,$$ 

\noindent 
then there exists $r'\in (0,r)$ such that $A+F$ is reducible in $C^\omega_{r'}(2\mathbb{T}^d )$.

\end{thm}

\dem Let $a\in [\max(1-\frac{1}{14^2},1-\frac{1}{G\cdot g(2)^2}),1[$. Let $ \epsilon_0>0,(\epsilon_n)_{n\in \mathbb{N}},
(N_n)_{n\in \mathbb{N}}$ as defined at the beginning of section \ref{iteration}.
Let $ (r_n),(\alpha_n),(m_n),(A_n),(F_n),(Z_n)$ be the sequences given by Proposition \ref{red0}. 

\noindent 
The sequence $(A_n)$ is bounded in 
$sl(2,\mathbb{R})$ for the operator norm so taking a subsequence $(A_{n_k})$, we find that $A_{n_k}$ tends to some $A_\infty\in gl(2,\mathbb{R})$. Now $\rho(A_\infty)$ is the limit of $\rho(A_{n_k})$ (see \cite{El92}, Lemma A.3) which implies that for all $n$,

$$\rho(A_\infty)=\rho(A_{n+1})-\lim_{k\rightarrow \infty}\sum_{j=n+1}^{n_k-1}(\rho(A_j)-\rho(A_{j+1})).$$

\noindent Moreover, 

\begin{equation}\rho(A+F)=\rho(A_\infty)+\pi\sum_{j\geq 0}\langle m_j,\omega\rangle
\end{equation}

\bigskip
\noindent 
(see also \cite{El92}). Therefore

\begin{equation}\begin{split}\mid \rho(A+F)-\pi \sum_{j\leq n}\langle m_j,\omega\rangle\mid
&=\mid  \rho(A_\infty)+\pi\sum_{j\geq n+1} \langle m_j,\omega\rangle\mid\\
&\leq \mid \alpha_{n+1}\mid +\sum_{j\geq n+1} \mid \alpha_j-\pi\langle m_j,\omega\rangle-\alpha_{j+1}\mid \\
&\leq \mid \alpha_{n+1}\mid +\sum_{j\geq n+1} \sqrt{\epsilon_j }
\end{split}\end{equation}

\noindent 
Suppose $\rho(A+F) $ satisfies


$$\forall m\in \mathbb{Z}^d, \ \mid \rho(A+F)-\pi \langle m,\omega\rangle\mid\geq \frac{\kappa'}{g(\mid m\mid)}.$$

\noindent 
In particular,


$$\mid {\rho}(A+F)-\pi  \sum_{j\leq n}\langle m_j,\omega\rangle\mid\geq \frac{\kappa'}{g(\mid  \sum_{j\leq n} m_j\mid) }$$

\noindent 
and so

\begin{equation} \frac{\kappa'}{g(\mid  \sum_{j\leq n} m_j\mid)}\leq \sum_{j\geq n+1} \sqrt{\epsilon_j}+\mid \alpha_{n+1}\mid.
\end{equation}

\noindent 
Let $n>n_0$. Assume $m_{n}\neq 0$. Then we have

$$\mid \alpha_{n+1}\mid\leq \mid \alpha_n - \pi \langle m_n,\omega\rangle\mid+\mid \alpha_n - \pi \langle m_n,\omega\rangle - \alpha_{n+1}\mid \leq \frac{\kappa}{4G(N_n)}+\sqrt{\epsilon_n}$$

\noindent so

\begin{equation}\label{croiss}  \kappa'\leq \left[\sum_{j\geq {n}}\sqrt{\epsilon_j}
+\frac{\kappa }{4G(N_{n}) }\right]g(\mid  \sum_{j\leq {n}} m_j\mid).
\end{equation}

\noindent 
Thus

\begin{equation}\kappa'\leq \left[\sum_{j\geq {n}}\sqrt{\epsilon_j}+\frac{\kappa }{4G(N_{n}) }\right] g(\sum_{j\leq {n}}\mid  m_j\mid) .
\end{equation}

\noindent 
Now

$$\sum_{j\geq {n}}\sqrt{\epsilon_j}=\frac{1}{1-(1-a)^{\frac{1}{4}}}\sqrt{\epsilon_{n}}\leq 2\sqrt{\epsilon_{n}}$$

\noindent 
and since, by definition of $N_n$,

$$
\epsilon_n
\leq \frac{(1-a)^2\kappa^2}{4 G(N_n)^2g(N_n)^2 }$$

\noindent 
then

\begin{equation}\kappa'\leq \frac{\kappa }{G(N_n) }g( \sum_{j\leq n}\mid  m_{j}\mid ).
\end{equation}

\noindent 
Now, note that $\sum_{j\leq n} \mid  m_{j}\mid \leq N_n^2$. This comes from the fact that, denoting by $m_{j_k} $ the subsequence of non-zero $m_j$'s, then for all $k$, 

\begin{equation}\begin{split}\mid \langle m_{j_{k+1}},\omega\rangle\mid &<
\mid \langle m_{j_{k+1}},\omega\rangle - \alpha_{j_{k+1}}\mid
+\mid \langle m_{j_k},\omega\rangle - \alpha_{j_k}+\alpha_{j_{k+1}}\mid
+\mid \langle m_{j_k},\omega\rangle - \alpha_{j_k}\mid\\
&< \frac{\kappa}{4G(N_{j_{k+1}})} + 2\sqrt{\epsilon_{j_k}}+\frac{\kappa}{4G(N_{j_k})}\\
& \leq \frac{\kappa}{4G(N_{j_{k+1}})} + \frac{\kappa}{2G(N_{j_k})}+\frac{\kappa}{4G(N_{j_k})}
\end{split}\end{equation}

\noindent 
which together with the arithmetic condition on $\omega$ implies that $N_{j_{k+1}}>N_{j_k}\geq k$. Therefore

\begin{equation}\kappa'
\leq \frac{\kappa }{G(N_n) }g(  N_n^2)
\end{equation}

\noindent 
and since 
by assumption $\frac{g( t^2)}{G(t)}$ is bounded,

\begin{equation}\kappa'\leq \kappa\sup_{t\geq n} \frac{g(  t^2)}{G(t)}.
\end{equation}

\noindent 
In other words, if

\begin{equation}\kappa'> \kappa \sup_{t\geq n} \frac{g(  t^2)}{G(t)}
\end{equation}

\noindent 
then $m_n=0$ for all $n\geq n_0$; and so $A+F$ is analytically reducible. $\Box$

\bigskip
\noindent 
This proves Theorem \ref{main}. 
Here is an easy consequence of the main result:

\begin{cor}\label{redsansn0} Let $A\in sl(2,\mathbb{R})$, $r>0$, $F\in C^\omega_r(\mathbb{T}^d )$. 
Assume

\begin{enumerate}


\item the map $t\mapsto \frac{g(t^2)}{G(t)}$ tends to $0$, 

\item $\rho(A+F)\in \mathrm{NR}_\omega(\kappa',g)$ for some $\kappa'>0$.

\end{enumerate}

\noindent 
There exist $\epsilon_0$ depending only on $g,\kappa,G,\rho(A+F), r$ such that if 

$$\mid F\mid_r\leq \epsilon_0,$$ 

%
\noindent 
then there exists $r'\in (0,r)$ such that $A+F$ is reducible in $C^\omega_{r'}(2\mathbb{T}^d )$.

\end{cor}

\dem Let $n_0$ be the smallest integer such that $\kappa'> \sup_{t\geq n_0}  \frac{g(t^2)}{G(t)}$. Take $\epsilon_0$ as in Theorem \ref{redu} so that it really depends on $g,\kappa,G,\rho(A+F), r$ and apply Corollary \ref{redu}. $\Box$

\bigskip
\noindent Theorem \ref{diophthm} is a particular case of Corollary \ref{redsansn0} since we can take $g(t)=t^\mu,G(t)=t^{\mu'}$ with $\mu'\geq \frac{\mu}{2}, \mu\geq 1,\mu'\geq 1$, as we will see in the next section.

\bigskip
\noindent Our main result will be made more convenient by the following corollary:

\begin{cor} Let $A\in sl(2,\mathbb{R}), r>0, F\in C^\omega_r(\mathbb{T}^d)$. Assume $\omega$ is a Brjuno vector and $\rho(A+F)$ is a $\frac{1}{2}$-Brjuno number with respect to $\omega$. There exists $\epsilon_0$ depending only on $\rho(A+F), \omega,r$ such that if $\mid F\mid_r\leq \epsilon_0$, then there exists $r'\in (0,r)$ such that $A+F$ is reducible in $C^\omega_{r'}$.
\end{cor}

\dem By assumption, there exists $\kappa'>0$ and $g$ positive increasing and continuous such that $\int_1^\infty \frac{\log g(t)}{t^{\frac{3}{2}}}dt<+\infty$ and $\rho(A+F)\in \mathrm{NR}_\omega(\kappa',g)$; there also exists $\kappa>0$ and $G'$ positive increasing and continuous such that $\int_1^\infty \frac{\log G'(t)}{t^2}dt<+\infty$ and $\omega\in\mathrm{NR}(\kappa,G')$. Now let for all $t, G(t)=t\max (G'(t), g(t^2))$. The function $G$ is positive increasing and continuous and $\omega\in \mathrm{NR}(\kappa,G)$. Since $\frac{g(t^2)}{G(t)}\leq \frac{1}{t}$ for all $t$, we can apply the previous corollary. $\Box$

\subsection{Possible choices of approximation functions}\label{Particular}

\noindent 
Here we give a few examples of approximation functions to which Theorem \ref{redu} can be applied.

\paragraph{Verification of Assumption \ref{rapportgG}}

\noindent 
Here are a few examples where Assumption \ref{rapportgG} holds, i.e $\frac{g(t^2)}{G(t)}$ is bounded:

\begin{enumerate}
\item \label{casdioph} $g(t)=t^\mu,G(t)=t^{\mu'}$ with $\mu'\geq \frac{\mu}{2}, \mu\geq 1,\mu'\geq 1$;
\item \label{casexp1} $g(t)=e^{t^\alpha}, G(t)=e^{t^{\alpha'}}$ with $\alpha\leq \frac{\alpha '}{2},\alpha<1,\alpha'<1$;
\item \label{casexp2} $g(t)=e^{t^\alpha},G(t)=e^{\frac{t}{(\log t)^\delta}},\alpha<1,\delta>1$.
\end{enumerate}

\noindent 
In the example \ref{casdioph}, and if $\mu'> \frac{\mu}{2}$, then, as noted in section \ref{Particular}, the condition on $\epsilon_0$ does not depend on $n_0$ and $\kappa'$ might be arbitrarily small, which corresponds to Eliasson's full-measure reducibility result in \cite{El92}.

\paragraph{Smallness conditions}

\noindent 
We shall make conditions \eqref{condepsilon} and \eqref{condepsilon2} more explicit for the particular cases that we mentioned before, namely, when $(g\cdot G)(t)=t^{\mu+\mu '}, \mu,\mu' >2$ (diophantine case), when $(g\cdot G)(t)=e^{t^\alpha+t^{\alpha'}}, \alpha,\alpha'<1$ and when $(g\cdot G)(t)=e^{\frac{t}{(\log t)^\delta}+t^\alpha},\delta>1,\alpha<1$.

%
%

\noindent 
Recall Condition \eqref{condepsilon}:

\begin{equation}\int_{(g\cdot G)^{-1}\left(\frac{\kappa}{2(1-a)^{\frac{n_0-5}{4}}\sqrt{\epsilon_0}}\right)}^\infty \frac{\log (g\cdot G)(t)}{t^2}dt\leq \frac{r_0}{4^{n_0+2}}
\tag{\ref{condepsilon}}\end{equation}

\begin{lem}If $(g\cdot G)(t)=t^{\mu+\mu'}$, Condition \eqref{condepsilon} is satisfied if

\begin{equation}(1-a)^{\frac{1}{8(\mu+\mu')}}\leq \frac{1}{2}
\  
\mathrm{and}\  
 \epsilon_0 \leq (\frac{r_0}{8(\mu+\mu')})^{4\mu}(1-a)^{\frac{3}{2}}\kappa
\end{equation}

\noindent 
or if

\begin{equation}\label{petitdioph} \epsilon_0 \leq (\frac{r_0}{4^{n_0+3}(\mu+\mu')})^{4(\mu+\mu')}\kappa.
\end{equation}

\end{lem}

\dem Rewrite Condition \eqref{condepsilon} as

\begin{equation}\int_{ b}^\infty \frac{(\mu+\mu')\log t}{t^2}dt\leq \frac{r_0}{4^{n_0+2}}
\end{equation}

\noindent 
where $b=\left(\frac{\kappa}{2(1-a)^{\frac{n_0-5}{4}}\sqrt{\epsilon_0}}\right)^{\frac{1}{(\mu+\mu')}} $. Integrating by parts, this is

\begin{equation}(\mu+\mu')\frac{\log b+1}{b}\leq \frac{r_0}{4^{n_0+2}}.
\end{equation}

\noindent 
It is enough that

\begin{equation}2(\mu+\mu')\frac{1}{\sqrt{b}}\leq \frac{r_0}{4^{n_0+2}}
\end{equation}

\noindent 
that is,

\begin{equation}\label{condbis}2(\mu+\mu') \left(\frac{2}{\kappa}(1-a)^{\frac{n_0-5}{4}}\sqrt{\epsilon_0}\right)^{\frac{1}{2(\mu+\mu')}}\leq \frac{r_0}{4^{n_0+2}}.
\end{equation}

\noindent 
which is true if \eqref{petitdioph} is satisfied. 
If moreover

\begin{equation}(1-a)^{\frac{1}{8(\mu+\mu')}}\leq \frac{1}{2}
\end{equation}

\noindent 
then \eqref{condbis} is satisfied as long as

\begin{equation} \epsilon_0 \leq (\frac{r_0}{8(\mu+\mu')})^{4(\mu+\mu')}(1-a)^{\frac{3}{2}}.\ \Box
\end{equation}


\begin{lem}If $(g\cdot G)(t)=e^{t^\alpha+t^{\alpha'}},\alpha'<\alpha<1$, then \eqref{condepsilon} holds if

\begin{equation} \label{talpha}  \epsilon_0  \leq \frac{\kappa}{4}\exp \left[-2\left(\frac{2\cdot 4^{n_0+2}}{r_0(1-\alpha)}\right)^{\frac{\alpha}{1-\alpha}}\right].
\end{equation}

\end{lem}

\dem By plugging $e^{t^\alpha+t^{\alpha'}}$ into \eqref{condepsilon} and recalling that $g\cdot G$ is increasing, one finds that 
%
%
%
%
%
%
%
%
%
%
%
%
Condition \eqref{condepsilon} holds if

\begin{equation}  (1-a)^{\frac{n_0-3}{2}}\epsilon_0  \leq \frac{\kappa}{4} \exp \left[-2\left(\frac{2\cdot 4^{n_0+2}}{r_0(1-\alpha)}\right)^{\frac{\alpha}{1-\alpha}}\right]
\end{equation}

\noindent 
so, in particular, \eqref{talpha} is a sufficient condition. $\ \Box$

\begin{lem}If $(g\cdot G)(t)=e^{\frac{t}{(\log t)^\delta}+t^\alpha}, \alpha<1, \delta>1$, then \eqref{condepsilon} holds if

\begin{equation}
\epsilon_0\leq \frac{\kappa}{4}\left((g\cdot G)\circ\exp\left[\left(\frac{4^{n_0+3}}{r_0(\delta-1)(1-\alpha) }\right)^{\frac{1}{(\delta-1)(1-\alpha)}}\right]\right)^{-2}.
\end{equation}
\end{lem}

\dem In this case, \eqref{condepsilon} can be rewritten

\begin{equation}\int_b^\infty \frac{1}{t^{2-\alpha}}dt+\int_{b}^\infty \frac{1}{t(\log t)^\delta}dt\leq \frac{r_0}{4^{n_0+2}}
\end{equation}

\noindent 
with $b=(g\cdot G)^{-1}\left(\frac{\kappa}{2(1-a)^{\frac{n_0-5}{4}}\sqrt{\epsilon_0}}\right)$. Integrating by parts, 
we compute

\begin{equation}\int_{b}^\infty \frac{1}{t(\log t)^\delta}dt=[ \frac{1}{(\log t)^{\delta-1}}]_b^\infty+\delta \int_{b}^\infty \frac{1}{t(\log t)^\delta}dt
\end{equation}

\noindent 
which implies

\begin{equation}(\delta-1)\int_{b}^\infty \frac{1}{t(\log t)^\delta}dt=\frac{1}{(\log b)^{\delta-1}}
\end{equation}

\noindent 
so that \eqref{condepsilon} is equivalent to

\begin{equation}
\frac{b^{\alpha-1}}{1-\alpha}+\frac{1}{ (\delta-1)(\log b)^{\delta-1}}\leq \frac{r_0}{4^{n_0+2}}.
\end{equation}

\noindent 
Now $\frac{1}{ (\delta-1)(\log b)^{\delta-1}}\leq \frac{r_0}{4^{n_0+3}}$ if

\begin{equation}
(g\cdot G)\circ\exp\left[\left(\frac{4^{n_0+3}}{r_0(\delta-1) }\right)^{\frac{1}{\delta-1}}\right]^2\leq \frac{\kappa}{2(1-a)^{\frac{n_0-5}{4}}\sqrt{\epsilon_0}}
\end{equation}

\noindent 
and $\frac{b^{\alpha-1}}{1-\alpha}\leq  \frac{r_0}{4^{n_0+3}}$ if

$$
\epsilon_0  \leq \frac{\kappa}{4}\exp \left[-2\left(\frac{ 4^{n_0+3}}{r_0(1-\alpha)}\right)^{\frac{\alpha}{1-\alpha}}\right]$$

\noindent 
so that \eqref{condepsilon} holds if 

\begin{equation}
\epsilon_0\leq \frac{\kappa}{4}\left((g\cdot G)\circ\exp\left[\left(\frac{4^{n_0+3}}{r_0(\delta-1)(1-\alpha) }\right)^{\frac{1}{(\delta-1)(1-\alpha)}}\right]\right)^{-2}.\ \Box
\end{equation}

\bigskip
\noindent Finally we consider Condition \eqref{condepsilon2}: first note that in these examples, one can drop the term $\sup_{t\in [1,n_0] } \frac{\log (G\cdot g)(t+1)}{t}$ in \eqref{condepsilon2} since it is a non increasing function. Therefore \eqref{condepsilon2} holds if, for instance,

\begin{equation}eC'\epsilon_0^{\frac{r_0}{4^{n_0+5}}}\leq (1-a)^2\kappa
\end{equation}

\noindent Thus in example \ref{casdioph}, \eqref{condepsilon2} holds for a suitable $a$ if 

\begin{equation}eC'\epsilon_0^{\frac{r_0}{4^{n_0+5}}}\leq  \frac{1}{2^{2(\mu+\mu')}}\kappa.
\end{equation}

\noindent
In example \ref{casexp1}, \eqref{condepsilon2} holds for a suitable $a$ if 

\begin{equation}eC'\epsilon_0^{\frac{r_0}{4^{n_0+5}}}\leq
\kappa \exp (-4 \cdot 4^{\frac{1}{1-\alpha'}})
\end{equation}

\noindent
and finally in example \ref{casexp2}, it holds under the analogous condition

\begin{equation}eC'\epsilon_0^{\frac{r_0}{4^{n_0+5}}}\leq
\kappa \exp (-4 \cdot 4^{\frac{1}{1-\alpha}}).
\end{equation}


\begin{thebibliography}{2}

\bibitem{AFK} Avila,A.; Fayad,B.; Krikorian,R.: \textit{A KAM scheme for $SL(2,R)$ cocycles with liouvillean frequency}, ArXiv, 1001.2878v1 (2010)

\bibitem{B} Bruno, A.D.: \textit{An analytic form of differential equations}, Mathematical Notes, 1969, 6:6, 927Ð931

\bibitem{C2}Chavaudret, C. \textit{Strong almost reducibility for analytic and Gevrey quasi-periodic cocycles}, arXiv:0912.4814 (2010)- to appear in Bulletin de la Socit Mathmatique de France

\bibitem{El92}Eliasson, L.H.:
\textit{Floquet solutions for the 1-dimensional quasi-periodic Schrdinger equation},
{Comm. Math. Phys.} {146},
{447-482} ({1992})

\bibitem{El01}Eliasson, L.H.: \textit{Almost reducibility of linear quasi-periodic systems}, Proc. Sympos. Pure Math. 69, 
679-705 (2001)

\bibitem{GM} Giorgilli, A.; Marmi, S.: \textit{Convergence radius in the Poincar-Siegel problem}, Discrete and Continuous Dynamical Systems, Series S, vol.3, n¡4, december 2010, 601-621

\bibitem{JM} Johnson,R.; Moser, J.: \textit{The rotation number for almost periodic potentials}, Comm. Math. Phys. 84 (1982), 403-438

\bibitem{MMY} Marmi,S.; Moussa,P.; Yoccoz,J.C.:\textit{The Brjuno functions and their regularity properties}, Comm. Math. Phys. 186 (1997), 265-293

\bibitem{P} Pschel, J. : \textit{KAM  la R}, arXiv: 0909.1015v2

\bibitem{R} Rssmann, H.: \textit{KAM iteration with nearly infinitely small steps in dynamical systems of polynomial character}, Discrete and Continuous Dynamical Systems, Series S, vol. 3, n¡4, december 2010, 683-718

\bibitem{Y} Yoccoz, J.C.: \textit{Petits diviseurs en dimension 1}, S.M.F., Astrisque 231(1995)

\bibitem{Y93} Young,L.S.: \textit{Lyapunov exponents for some quasi-periodic cocycles}, Ergodic Theory and Dynamical Systems 17 (1997) n¡2, 483-504

\bibitem{ZW} Wang,J.; Zhou,Q.: \textit{Reducibility results for quasiperiodic cocycles with liouvillean frequency}, J. Dynam. Differential Equations, DOI: 10.1007/s10884-011-9235-0

\end{thebibliography}
\end{document}